# Revisiting the Algebraic and Analytic Descriptions of Quantum Mechanics

**Ortwin Fromm, Felicitas Ehlen**

## ABSTRACT

The Wintner–Wielandt theorem exposes a fundamental incompatibility between the canonical commutation relation (CCR) and the Hilbert-space formulation of quantum mechanics. While the conventional approach preserves the CCR by extending the Hilbert-space framework to a rigged Hilbert space, this work develops the complementary alternative. Retaining the Hilbert space, the unavoidable trace obstruction requires a boundary modification of the CCR, consistently realizable only in a space of finite, albeit arbitrarily large, dimension. In this setting, the eigenvalue problems of position and momentum admit explicit solutions within the Hilbert space. We show that the highest-energy components of the position eigenvectors retain a non-vanishing uniform weight, implying that the full spectrum remains relevant. The resulting autonomous matrix algebra also yields a new physical observable, the *action difference*. In the position basis, this observable gives rise to a momentum matrix with a discrete Cauchy kernel. After an alternating phase transformation, this matrix acts on plane waves in the infinite equidistant lattice model as the Schrödinger momentum operator, revealing a structural connection between the algebraic and analytic approaches. This framework retains the standard Hilbert-space spectral theory and remains closely connected to the familiar analytical structure of quantum mechanics, including the momentum operator and the uncertainty principle.

## 1. INTRODUCTION

In 1925, Werner Heisenberg, Max Born, and Pascual Jordan [1] devised matrix mechanics. Their fundamental formalism was based on tridiagonal matrix operators $Q$ and $P$ for position and momentum, as well as the canonical commutation relation (CCR) $[Q,P]=i\cdot I$. However, it quickly became apparent that this relation is unrealizable in finite-dimensional spaces due to the inherent "trace defect." In 1926, Erwin Schrödinger [2], [3] developed wave mechanics. He assigned the multiplication operator $q$ to the position matrix $Q$ and the differential operator $-i\cdot\partial/\partial q$ to the momentum matrix $P$. This assignment allowed the validity of the CCR to be realized consistently on a space of smooth functions in the continuum for the first time. John von Neumann [4] brought mathematical rigor to the formalism in 1932. He demonstrated that both the space $\ell^2$ of square-summable sequences (matrix mechanics) and the space $L^2(\mathbb{R})$ of square-integrable functions (wave mechanics) are mutually isomorphic, complete Hilbert spaces, while the Stone-von Neumann theorem [5], [6] establishes the unitary equivalence of all regular irreducible representations of the canonical commutation relation. Only this Hilbert space theory enabled a consistent statistical interpretation. At this historical juncture, the foundation of quantum mechanics thus rested on two pillars: the CCR as the central algebraic structural relation and the Hilbert space geometry as the foundation for the probabilistic interpretation.

The Wintner-Wielandt theorem [7], [8], however, revealed a fundamental structural tension between these two principles. It proved that the classical CCR is incompatible with the Hilbert space structure in the infinite-dimensional setting: Since the theorem excludes bounded operators for the validity of the CCR on a Hilbert space, the operators $Q$ and $P$ must be unbounded. As unbounded operators, however, they are not defined on the entire Hilbert space $L^2(\mathbb{R})$. A common invariant domain is provided by the Schwartz space $S(\mathbb{R})$. The unrestricted validity of the CCR is therefore incompatible with the purely Hilbert-space-based formulation. Consequently, the Wintner-Wielandt theorem naturally leads to two complementary scenarios for the description of physical reality:

Either one chooses the unmodified CCR and relinquishes the pure Hilbert space framework, or one preserves the Hilbert space framework and modifies the CCR.

In the first scenario, the CCR is taken as the fundamental axiom. The classical route proceeds via the Schrödinger representation to the analytic setting, where the commutation relation holds exactly on the Schwartz space $S(\mathbb{R})$, regarded as a dense subspace of $L^2(\mathbb{R})$. Since $S(\mathbb{R})$ itself is not a Hilbert space and is insufficient for a complete spectral treatment of the unbounded operators $Q$ and $P$, one is naturally led to the Gelfand triple $S(\mathbb{R}) \subset L^2(\mathbb{R}) \subset S'(\mathbb{R})$. By admitting generalized eigenvectors as tempered distributions in $S'(\mathbb{R})$, the pure Hilbert-space framework is therefore deliberately replaced by an extended formalism.

In the second scenario, by contrast, the preservation of a pure Hilbert-space framework is taken as the primary requirement. This naturally leads to the sequence space $\ell^2$. Although the rigged Hilbert-space triple $s \subset l^2 \subset s'$, obtained as the unitary image of $S \subset L^2 \subset S'$ also provides a realization of the CCR, its introduction is inherited from the analytical formalism rather than required by the algebraic one. In the present work, we therefore retain the Hilbert-space framework and modify the CCR instead.

Since the tridiagonal structure of the operators $Q$ and $P$ shows that any modification of the CCR affects only the boundary term, such a modification cannot be realized in the boundaryless space $\ell^2$. This motivates the transition to the finite-dimensional Hilbert space, where the boundary is well defined. By indexing the basis vectors as $n = 0,1,...,N$, we arrive at the Hilbert space $\mathbb{C}^{N+1}$ with a finite but arbitrarily large dimension. Within this space, all algebraic operations are well defined. Furthermore, all operators are bounded and representable as finite-dimensional matrices.

The modified CCR naturally gives rise to a new observable, $D_N = [Q_N, P_N]/i$ which we term the *action difference*. This quantity was first introduced in our earlier work [9]. The tridiagonal structure of the position and momentum matrices yields the diagonal form $D_N = diag(1,1,...,-N)$, which immediately manifests the vanishing of the trace. Notably, the boundary term $-N$ scales unboundedly with the dimension of the space, suggesting that the action difference may represent a genuine physical observable. The modified CCR addresses a long-standing challenge in quantum mechanics: within a finite-dimensional space, the commutator cannot, in principle, coincide exactly with the identity operator. This resulting structural constraint and its physical implications continue to be actively discussed in modern mathematical physics, as demonstrated by the works of de la Torre and Goyeneche [10] as well as Santhanam [11].

The principal advantage of the algebraic approach lies in the unrestricted applicability of all matrix operations within the Hilbert space $\mathbb{C}^{N+1}$. In particular, the eigenvalue problems for the position and momentum operators admit exact finite-dimensional solutions. We show that the spectra of $Q$ and $P$ lie densely on the real axis as $N$ increases. Furthermore, we prove that the squares of the final components of all position eigenvectors possess the identical value of $1/(N+1)$, which corresponds exactly to the mean of the squares of all remaining components. Consequently, the energetically highest components retain their relative weight as $N$ grows. Within the present framework, this implies that during a sharp position measurement, the entire energetic spectrum up to the highest components is physically relevant. Thus the finite measurement window of a finite observer becomes vanishingly small compared with the arbitrarily large space $\mathbb{C}^{N+1}$. Almost all information is therefore inaccessible to the observer. This restricted accessibility provides a structural motivation for interpreting quantum randomness in epistemic terms, in close analogy with the statistical descriptions of Boltzmann and Gibbs. The present framework therefore naturally supports an epistemic interpretation of quantum randomness. The epistemic limitation may be understood as an intrinsic form of unavoidable decoherence arising at the interface between the finite observational window and the arbitrarily large space. Despite their different mathematical formulations, the finite algebraic and continuous analytic approaches share essential structural features. In particular, the standard Hilbert-space spectral theory remains fully applicable, while the momentum operator and the uncertainty principle retain their familiar significance within the modified algebraic framework. Furthermore, we demonstrate that, when transformed into the position basis, the action difference assumes the dimension-independent form $D^{(q)} = I - J$, where $I$ denotes the identity operator and $J$ the all-ones operator. We term $D^{(q)}$ the *Leibniz operator*, as it couples every spatial state equally to every other state. This invariant structure gives rise to a matrix representation of the momentum operator with a discrete Cauchy kernel, providing the algebraic counterpart of the differential momentum operator in the Schrödinger coordinate representation. As demonstrated in Appendix C, an alternating phase transformation of this Cauchy matrix reproduces the action of the Schrödinger momentum operator on plane waves in the infinite equidistant lattice model exactly.

This work deepens and extends the foundations published under the same title in arXiv:2602.14126.

With the choice of the matrix-algebraic path, the subsequent development in this work follows naturally:

In Section 2, we systematically develop the formalism and present our main results:

Starting from the matrix representations of position and momentum (2.1) and the harmonic oscillator (2.2) for exemplification, we derive formulas for the eigenvalues of $Q_N$ and $P_N$ (2.3) and evaluate their spectra (2.4). Subsequently, we deduce the formulas for the eigenvector components (2.5) as well as their respective final components (2.6). We then confirm the unitary equivalence of $Q_N$ and $P_N$ (2.7). On this foundation, we formulate the matrix of the action difference $D_N$ (2.8), which can then be transformed into the position basis (2.9), giving rise to the Leibniz operator. Utilizing the action difference, we compute the position representation of momentum and obtain the discrete Cauchy kernel (2.10).

Finally, Section 3 offers comprehensive discussions and algebraic reflections, contrasting this finite approach with the continuous formalism, discussing the epistemic interpretation of quantum randomness, and investigating the emergence of relational space. Appendix C further examines the connection to the Schrödinger representation through the momentum operator in the position basis.

**Remark**: In conventional truncated formulations, $N+1$ serves as the numerical cutoff dimension of an underlying infinite-dimensional Hilbert space. In contrast, no such truncation is introduced along the matrix-algebraic path pursued here. Motivated by the Wintner-Wielandt theorem as the second alternative, $N+1$ labels the arbitrarily large dimension of the finite-dimensional Hilbert-space representation.

## 2. FORMALISM AND RESULTS

### 2.1 Matrix Representation of Position and Momentum

In their groundbreaking 1925 paper *„Zur Quantenmechanik"* [1], Max Born and Pascual Jordan introduced infinite-dimensional tridiagonal matrix representations for the position and the momentum operators $Q$ and $P$. Within the autonomous finite-dimensional Hilbert space $\mathbb{C}^{N+1}$ (spanned by the basis states $e_l$, $l = 0,1,...,N$), these fundamental operators take an explicit bounded form. To simplify the algebraic exposition, we adopt natural units throughout this work, setting the mass, the oscillator frequency, and the reduced Planck constant to unity: $m = \omega = \hbar = 1$. Consequently, all position, momentum, and energy variables are rendered dimensionless. In the standard Fock basis, the operators are defined by their matrix elements, which we explicitly display below for the case of $N = 3$, corresponding to dimension 4:

$$Q_3 = \frac{1}{\sqrt{2}}\begin{pmatrix} 0 & \sqrt{1} & 0 & 0 \\ \sqrt{1} & 0 & \sqrt{2} & 0 \\ 0 & \sqrt{2} & 0 & \sqrt{3} \\ 0 & 0 & \sqrt{3} & 0 \end{pmatrix}, \quad P_3 = \frac{i}{\sqrt{2}}\begin{pmatrix} 0 & -\sqrt{1} & 0 & 0 \\ \sqrt{1} & 0 & -\sqrt{2} & 0 \\ 0 & \sqrt{2} & 0 & -\sqrt{3} \\ 0 & 0 & \sqrt{3} & 0 \end{pmatrix} \tag{1}$$

### 2.2 Hamiltonian of the Harmonic Oscillator and Energy Eigenvalues

The matrix of the Hamiltonian $H_N = \left(Q_N^2 + P_N^2\right)/2$ is diagonal. The only deviation from equally spaced harmonic-oscillator spectrum occurs at the terminal entry (index $N$), whose value is shifted from $(2N+1)/2$ to $N/2$: $H_N = diag\left(1/2,\ 3/2,\ ...,(2N-1)/2,\ N/2\right)$.

For the case $N = 3$ (dimension 4) we have $H_3 = \frac{1}{2}\left(Q_3^2 + P_3^2\right) = \begin{pmatrix} 1/2 & 0 & 0 & 0 \\ 0 & 3/2 & 0 & 0 \\ 0 & 0 & 5/2 & 0 \\ 0 & 0 & 0 & 3/2 \end{pmatrix}$.

Thus, the terminal eigenvalue lies well inside the spectral band and contributes only a fraction $1/\left(N+1\right)$ of the normalized spectrum, which vanishes as $N \to \infty$.

### 2.3 Eigenvalues of *Q* and *P*

Expanding the characteristic polynomial $\det(Q_N - \lambda I_N)$ along rows or columns yields a recursion identical to that of the Hermite polynomials $H_n(x)$. Comparison with the Hermite recurrence shows that the characteristic polynomial of $Q_N$ is proportional to $H_{N+1}(\lambda)$. Hence its zeros are precisely the eigenvalues $\lambda_0, \lambda_1, ..., \lambda_N$ of the position matrix $Q_N$ (cf. Szegő [12] or Abramowitz and Stegun [13]). The characteristic polynomial of the momentum matrix $P_N$ satisfies the same recurrence relation as that of the position matrix $Q_N$. Accordingly, we use the same notation $\lambda_l$ for both spectra

$$\lambda_l = x_l = p_l\,,\ l = 0,..,N\ , \tag{2}$$

meaning that the eigenvalues of $Q_N$ and $P_N$ are identical and given by the zeros of the Hermite polynomials $H_{N+1}(\lambda)$, establishing that $Q_N$ and $P_N$ are unitarily equivalent. The well-known unitary Fourier equivalence between $Q$ and $P$ in the infinite-dimensional case thus has its counterpart in the finite-dimensional setting.

### 2.4 Dense Spectra

As $N \to \infty$, classical results of Stieltjes on the zeros of Hermite polynomials [12] imply that the eigenvalues of $Q_N$ and $P_N$ become asymptotically dense on the real line. Consequently, the finite-dimensional representation in the autonomous Hilbert space $\mathbb{C}^{N+1}$ reproduces the spectral distribution of the continuous Schrödinger representation in the large-$N$ limit, without requiring an explicit embedding into an infinite-dimensional Hilbert space.

### 2.5 Closed-Form Expression for the Eigenvector Components

Let $\overline{u}_l^{(N)}$ denote an unnormalized eigenvector of $Q_N$ (see 2.1) corresponding to the eigenvalue $\lambda_l$; $l = 0,1,...,N$. Setting $\overline{u}_{l,0}^{(N)} = 1$ we obtain the following recurrence relations, from which all components are obtained recursively.

$$\sqrt{1}\cdot\overline{u}_{l,1}^{(N)} = \sqrt{2}\cdot\lambda_l\cdot\overline{u}_{l,0}^{(N)}$$
$$\sqrt{2}\cdot\overline{u}_{l,2}^{(N)} = \sqrt{2}\cdot\lambda_l\cdot\overline{u}_{l,1}^{(N)} - \sqrt{1}\cdot\overline{u}_{l,0}^{(N)}$$
$$\sqrt{3}\cdot\overline{u}_{l,3}^{(N)} = \sqrt{2}\cdot\lambda_l\cdot\overline{u}_{l,2}^{(N)} - \sqrt{2}\cdot\overline{u}_{l,1}^{(N)}$$
$$\ldots$$
$$\sqrt{k+1}\cdot\overline{u}_{l,k+1}^{(N)} = \sqrt{2}\cdot\lambda_l\cdot\overline{u}_{l,k}^{(N)} - \sqrt{k}\cdot\overline{u}_{l,k-1}^{(N)}$$

In the Appendix A, the result $\overline{u}_{l,k}^{(N)}(\lambda_l) = H_k(\lambda_l)/\sqrt{2^k\cdot k!}$ for the components is derived from the recurrence relations via mathematical induction, where $H_k(\lambda)$ are again the Hermite polynomials ( $H_0(\lambda)=1,\ H_1(\lambda)=2\lambda,\ H_2(\lambda)=4\lambda^2-2,...$ ). Normalization yields the eigenvector $u_0^{(N)}, u_1^{(N)},..., u_N^{(N)}$ with the components

$$u_{l,k}^{(N)}(\lambda_l) = \frac{H_k(\lambda_l)}{\sqrt{2^k\cdot k!}} \Big/ \sqrt{\sum_{j=0}^{N}\frac{(H_j(\lambda_l))^2}{2^j\cdot j!}}\,. \tag{3}$$

### 2.6 Constant Final Components

The individual components of the normalized position eigenvectors $u_l^{(N)}$, $l = 0,1,...,N$ depend on the eigenvalues $\lambda_l$, which are generally obtained numerically. The algebraic boundary, however, is governed by an exact symmetry, stated here as the *Invariance Theorem*: the final component of each normalized position eigenvector is strictly invariant with respect to the position index $l$ and takes the constant value:

$$u_{l,N}^{(N)} = \frac{1}{\sqrt{N+1}}, \text{ for all } l = 0,..,N\,. \tag{4}$$

**Proof:** From the explicit expression for the normalized components derived in Section 2.5 and proved by induction in Appendix A, the terminal boundary element is determined by evaluating the Hermite polynomials at the roots $\lambda_l$. The resulting summation in the denominator can be evaluated analytically via the classical Christoffel-Darboux identity [13], as detailed in Appendix B: $\sum_{j=0}^{N}\frac{H_j(\lambda_l)^2}{2^j\cdot j!}=\frac{(N+1)\cdot H_N(\lambda_l)^2}{2^N\cdot N!}$.

By substituting this identity into the normalization condition, the explicit dependence on the eigenvalue $\lambda_l$ completely cancels out for every root:

Thus all *N*-components are $u_{l,N}^{(N)}=\frac{\frac{H_N(\lambda_l)}{\sqrt{2^N\cdot N!}}}{\sqrt{\sum_{j=0}^{N}\frac{H_j(\lambda_j)^2}{2^j\cdot j!}}}=\frac{\frac{H_N(\lambda_l)}{\sqrt{2^N\cdot N!}}}{\sqrt{\frac{(N+1)\cdot H_N(\lambda_l)^2}{2^N\cdot N!}}}=\frac{1}{\sqrt{N+1}}$.

This proves the invariance of the terminal components across the entire position spectrum. Consequently, the modal matrix is given by $U_N$ with strictly constant values in the final row.

$$U_N=\left[u_0^{(N)},u_1^{(N)},\dots,u_N^{(N)}\right]. \tag{5}$$

**Remark**: Alternatively, the invariance theorem $1/\sqrt{N+1}$ was originally derived [9] via a geometric approach using a Koecher rotation [14]. In this framework, the corresponding modal matrix $U$ is partitioned into blocks, shifting the proof from an algebraic verification to a geometric representation.

For clarity, we will display the modal matrix $U$ for the case $N=3$ (dimension 4)

$$U=\begin{pmatrix}-0.214 & 0.674 & -0.674 & 0.214\\ 0.5 & -0.5 & -0.5 & 0.5\\ -0.674 & -0.214 & 0.214 & 0.674\\ 0.5 & 0.5 & 0.5 & 0.5\end{pmatrix}.$$

**Consequence:** Only the terminal row of the modal matrix $U_N$ in Eq. (5) is needed to transform the action difference into the position basis. This invariant constant last row simplifies the boundary algebra and is essential for the analysis in section 2.9.

## 2.7 Unitary Equivalence between *Q* and *P*

As demonstrated in Section 2.3, the characteristic polynomials of both the position matrix $Q_N$ and the momentum matrix $P_N$ satisfy the identical three-term recurrence relation of the Hermite polynomials $H_{N+1}(\lambda)$. Consequently, their discrete spectra are identical: $\sigma(Q_N)=\sigma(P_N)$. According to the spectral theorem for finite-dimensional normal operators, two self-adjoint matrices sharing the same set of distinct eigenvalues are unitarily equivalent: $U_N^\dagger Q_N U_N=V_N^\dagger P_N V_N$. Therefore, there exists a finite-dimensional unitary transformation matrix $F_N=U_N V_N^\dagger$ that maps the position representation directly onto the momentum representation, $P_N = F_N^\dagger Q_N F_N$ . The operator $F_N$ serves as the exact finite-dimensional counterpart to the continuous quantum mechanical Fourier transform.

### 2.8 Action Difference

Our algebraic framework is formulated in terms of finite-dimensional matrices $Q_N$ and $P_N$, where $N \in \mathbb{N}$ can be chosen arbitrarily large. Throughout this work, the subscript $N$ indicates that the operators act on an $(N+1)$-dimensional space, spanned by the Fock states $l = 0,1,...,N$ . Following Born & Jordan's original foundational construction, these matrices are tridiagonal (see Section 2.1). Their commutator defines a Hermitian operator:

$$D_N = [Q_N, P_N]/i, \tag{6}$$

which we denote as the *action difference*. This operator extends the canonical pair by a third observable with the dimension of action, which acquires a non-trivial form only within the finite algebraic framework. In the Fock basis, this operator takes the explicit diagonal form:

$$D_N = diag(1,1,...,-N), \tag{7}$$

thereby exhibiting an intrinsic asymmetry localized entirely in its final diagonal entry.
Having previously introduced $Q_3$ and $P_3$ in Section 2.1, we explicitly display the corresponding action difference $D_3$ in its full matrix representation:

$$D_3 = \frac{Q_3 \cdot P_3 - P_3 \cdot Q_3}{i} = diag(1,1,1,-3) = \begin{pmatrix} 1 & 0 & 0 & 0 \\ 0 & 1 & 0 & 0 \\ 0 & 0 & 1 & 0 \\ 0 & 0 & 0 & -3 \end{pmatrix}.$$

From the general structure in Eq. (7), it follows that the operator can be decomposed as:

$$D_N = I_N - (N+1) \cdot e_N \cdot e_n^\dagger, \tag{8}$$

where the second term represents the boundary term enforced by the algebraic closure within the finite framework. Here, $I_N$ denotes the $(N+1)$-dimensional identity matrix, and $e_N$ is the $(N+1)$-dimensional unit vector directed along the terminal coordinate axis (index $N$ ).
The action difference also enters directly into the Heisenberg uncertainty relation. For any normalized state $\psi \in \mathbb{C}^{N+1}$ the modified commutation relation $[Q_N, P_N] = iD_N$ implies $\Delta_\psi Q_N \, \Delta_\psi P_N \geq \frac{1}{2} \left| \langle \psi, D_N \psi \rangle \right|$ . Using the decomposition in Eq. (8), this yields $\Delta_\psi Q_N \, \Delta_\psi P_N \geq \frac{1}{2} \left| 1 - (N+1) |c_N|^2 \right|$, where $c_N = \langle e_N, \psi \rangle$ denotes the final component of the normalized state in the Fock basis. This inequality holds exactly for every finite dimension $N+1$ and shows that the action difference determines the state-dependent lower bound of the uncertainty product.

## 2.9 The Leibniz Operator and Relational Space

We now perform a standard basis transformation to express the action difference $D$ in the position representation:

$$D_N^{(q)} = U_N^\dagger D_N U_N \ , \qquad (9)$$

where $U_N$ (Eq. (5)) denotes the unitary modal matrix formed by the eigenvectors of $Q_N$. Utilizing the algebraic decomposition established in Eq. (8) of Section 2.8, we can transform the two structural terms independently:

$$D_N^{(q)} = U_N^\dagger \left( I_N - (N+1) \cdot e_N e_N^T \right) U_N = I_N - (N+1) \cdot U_N^\dagger \cdot e_N e_N^T \cdot U_N \ .$$

We note that the dyadic product $e_N e_N^T$ isolates the final diagonal entry of the matrix in the Fock basis. Consequently, the transformed matrix product $(N+1) \cdot U_N^\dagger e_N e_N^\dagger U_N$ is composed entirely of the boundary components $u_{l,N}^{(N)}$ of the normalized eigenvectors.

At this critical juncture, our fundamental invariance theorem from Eq. (4) in Section 2.6 applies: since $u_{l,N}^{(N)} = 1/\sqrt{N+1}$ holds universally for all spectral indices $l = 0,...,N$, every single entry of this transformed dyadic matrix evaluates precisely to the constant value $1/\sqrt{N+1}$. This constant value cancels out perfectly with the algebraic prefactor $(N+1)$. Therefore, the action difference transformed into the position basis simplifies to the forms:

$$D_N^{(q)} = I_N - J_N \qquad (10a) \qquad \text{and, in index-free notation} \qquad D^{(q)} = I - J \ , \qquad (10b)$$

where $J_N$ denotes the $(N+1)\times(N+1)$ all-ones matrix.

Crucially, the explicit dependence on the matrix dimension $N$ disappears completely from the right-hand side of Eq. (10a). This universal structure justifies the index-free notation $D^{(q)} = I - J$ which we shall refer to as the *Leibniz operator*, reflecting its intrinsically relational character. For clarity, we display the resulting matrix for the case $N = 3$, (dimension 4):

$$D_3^{(q)} = \begin{pmatrix} 0 & -1 & -1 & -1 \\ -1 & 0 & -1 & -1 \\ -1 & -1 & 0 & -1 \\ -1 & -1 & -1 & 0 \end{pmatrix}$$

Unlike the continuous Schrödinger representation, where the canonical commutation relation is represented by the identity operator in the position basis, the finite-dimensional framework transforms the boundary contribution into a matrix with vanishing diagonal and uniform off-diagonal entries. With respect to the action difference, no position state retains an independent algebraic status; instead, each is characterized solely through its relation to all others. the Leibniz operator $D^{(q)} = I - J$ thus provides a natural algebraic representation of the action difference in finite quantum space. Its universal form suggests an intrinsically relational structure, in which localisation is replaced by homogeneous relations between all position states.

## 2.10 Position Representation of the Momentum Matrix: The Discrete Cauchy Kernel

Building upon the invariant structure of the Leibniz operator $D^{(q)} = I - J$ derived in Section

2.9, we now evaluate the explicit matrix elements of the momentum operator in the position representation. This purely algebraic derivation yields a matrix with a discrete Cauchy kernel, providing the algebraic counterpart of the differential momentum operator in the Schrödinger coordinate representation and illustrating the structural consistency of the finite-dimensional framework.

Let $x_0, x_1, ..., x_N$ denote the $N+1$ distinct eigenvalues of the position matrix $Q_N$, such that its diagonalized representation is given by $\Lambda_N = diag(x_0, x_1, ..., x_N)$. Transforming the foundational commutator relation $i\, D_N = [Q_N, P_N]$ (Eq. 6) into the position basis leads to the following matrix identity:

$$iD_N^{(q)} = U_N^\dagger [Q_N, P_N] U_N = U_N^\dagger Q_N U_N \cdot U_N^\dagger P_N U_N - U_N^\dagger P_N U_N \cdot U_N^\dagger Q_N U_N = \Lambda_N \cdot P_N^{(q)} - P_N^{(q)} \cdot \Lambda_N \quad (10)$$

where $P_N^{(q)} = U_N^\dagger P_N U_N$ represents the momentum matrix in the position representation. Evaluating Eq. (10) elementwise yields:

$$\left(iD_N^{(q)}\right)_{j,k} = x_j \cdot P_{N;\,j,k}^{(q)} - P_{N;\,j,k}^{(q)} \cdot x_k\,, \quad (11)$$

respectively

$$P_{N;\,j,k}^{(q)} = \frac{\left(iD_N^{(q)}\right)_{j,k}}{x_j - x_k}. \quad (12)$$

At this point, the essential difference from the conventional infinite-dimensional Hilbert space formalism becomes apparent. If one were to adopt the standard canonical commutation relation $[Q,P] = i \cdot I$, the transformed action difference $D_N^{(q)}$ would incorrectly reduce to the identity matrix $I_N$. For the matrix elements this substitution forces the right-hand side of Eq. (12) to vanish: $P_{N;\,j,k}^{(q)} = 0$ for $j \neq k$, and for the diagonal element Eq. (11) yields the contradiction $i = 0$, which ultimately yields an algebraically inconsistent result.

In contrast, by preserving the intrinsic boundary term embedded within our finite-dimensional operator $D_N^{(q)} = I_N - J_N$, we have $D_{N;\,j,k}^{(q)} = -1$ for $j \neq k$ and $D_{N;\,j,j}^{(q)} = 0$. So Eq. (12) yields:

$$P_{N;\,j,k}^{(q)} = \frac{i}{x_k - x_j} \text{ for } j \neq k\,. \quad (13)$$

For the diagonal elements $(j = k)$ Eq. (11) yields the structurally indeterminate form $0 = 0$. From Eq. (1) we see that the matrix $P_N$ is antisymmetric: $P_N^T = -P_N$ and $Q_N$ is real symmetric. So the modal matrix $U_N$ of $Q_N$, satisfying $U_N^T = U_N^{-1}$, is real. Transforming $P_N$ into the position basis yields $P_N^{(q)} = U_N^T P_N U_N$ and $\left(P_N^{(q)}\right)^T = \left(U_N^T P_N U_N\right)^T = U_N^T P_N^T U_N =$ $= -U_N^T P_N U_N = -P_N^{(q)}$. Since $P_N^{(q)}$ is antisymmetric, all its diagonal entries vanish:

$$P_{N;\,j,j}^{(q)} = 0 \;\forall\, j = 0, ..., N \quad (14)$$

Equations (12) and (13) demonstrate that, in the algebraic position representation, the momentum operator is represented by a discrete Cauchy kernel. This kernel is the finite-dimensional analogue of the differential momentum operator in the continuous Schrödinger representation and naturally reflects the non-local algebraic structure of the finite-dimensional framework.

## 3. DISCUSSION AND ALGEBRAIC REFLECTIONS

### 3.1 Comparison of the Finite and Continuous Formalism

The algebraic framework developed in Section 2 establishes a self-contained description of quantum mechanics within a strictly finite-dimensional Hilbert space $\mathbb{C}^{N+1}$. While continuous Schrödinger wave mechanics necessarily alternates between restrictions to the Schwartz space $S(\mathbb{R})$ due to incompatible operator domains of $Q$ and $P$, and extensions to the space of tempered distributions $S'(\mathbb{R})$ outside $L^2(\mathbb{R})$, the formalism presented here remains strictly within the Hilbert space $\mathbb{C}^{N+1}$ throughout all operations.

The determination of the eigenvalues of $Q_N$ and $P_N$ via the common roots of the Hermite polynomials (Section 2.3) ensures a dense spectrum in the real line by applying Stieltjes' theorems (Section 2.4). From the identity of the characteristic polynomial, the strict unitary equivalence of these operators follows (Section 2.7). In contrast, in the continuous space, the unitary equivalence of position states (delta distributions) and momentum states (plane waves) is realized by the continuous Fourier integral, which transforms both into objects outside $L^2(\mathbb{R})$. The relation between the finite algebraic and continuous analytical formulations can be understood at two complementary levels. The finite-dimensional spaces $\mathbb{C}^{N+1}$ form an increasing sequence of subspaces whose union is dense in $l^2$, while $l^2$ and $L^2(\mathbb{R})$ are unitarily isomorphic as separable infinite-dimensional Hilbert spaces. This provides the abstract Hilbert-space connection between the two formulations.

The corresponding connection at the operator level is discussed in Appendix C. In the position representation, the position matrix takes the form of a multiplication operator, while the Leibniz form of the action difference gives rise to the Cauchy structure of the momentum matrix. After an alternating phase transformation, this structure reproduces the action of the Schrödinger momentum operator on plane waves in the idealized infinite equidistant lattice model. This reveals a further, more concrete connection between the algebraic and analytic formulations of quantum mechanics at the operator level.

The Heisenberg uncertainty relation provides a further connection between the two formulations. As derived in Section 2.8, the finite-dimensional algebraic framework yields the exact state-dependent inequality $\Delta_\psi Q_N \, \Delta_\psi P_N \geq \frac{1}{2}\left|1-(N+1)\left|c_N\right|^2\right|$. For states with finite Fock expansions, the conventional lower bound $1/2$ holds exactly once the highest occupied energy level lies below the algebraic boundary, while for finite-dimensional projections of Schwartz-class states it is recovered asymptotically. By contrast, every normalized position eigenstate satisfies $\left|c_N\right|^2 = 1/(N+1)$, as established in Section 2.6. Consequently, the expectation value of the action difference vanishes, and the lower bound becomes zero. This is also evident because the position uncertainty of an exact position eigenstate vanishes, yielding a zero uncertainty product at every finite dimension, while the momentum uncertainty remains finite but grows proportionally to $\sqrt{N}$ in the equidistant lattice approximation of Appendix C. Such a divergence parallels the increasing momentum uncertainty associated with sharp localization in Schrödinger wave mechanics. However, the corresponding exact position eigenstates in the Schrödinger representation lie outside $L^2(\mathbb{R})$, where the uncertainty product is undefined.

The decisive structural transition occurs in Sections 2.5 and 2.6: the recurrence formula for the components $u_{l,k}$ of the position eigenvectors $u_l$, derived by mathematical induction in Appendix A, can be evaluated for the highest index $k = N$ exactly by using the Christoffel-Darboux identity in Appendix B. The result demonstrates that the highest-energy components of all position eigenvectors remain stable relative to the remaining components; namely, the last component squared, $\left|u_{l,N}\right|^2 = 1/(N+1)$, is precisely the mathematical mean of all component squares over the vector. Thus, although the squared boundary components tend to zero as $N$ increases, they are not suppressed relative to the other components. This asymptotic behaviour of the finite boundary term stands in sharp mathematical contrast to the convergence properties of the infinite-dimensional function space when a sharp position measurement or a localized state preparation is modeled. In the continuous formalism, an exactly localized position state is represented by the generalized eigenvector $\delta(x-x_0)$, which does not belong to $L^2(\mathbb{R})$. Its expansion coefficients $c_n$ in the Fock basis decay asymptotically only as $n^{-1/4}$ for large quantum numbers, according to the Stirling approximation for the harmonic oscillator [13]. Consequently, the corresponding series $\sum_n \left|c_n\right|^2$ diverges, reflecting the non-normalizability of the idealized position eigenstate rather than that of a square-integrable quantum state. If, instead, the singular distribution is replaced by a square-integrable Gaussian wave packet with finite standard deviation $\sigma$, analyticity enforces a qualitatively different behaviour: the coefficients decay exponentially, proportional to $\exp(-2\sigma^{2n})$ [13]. The finite-dimensional matrix formalism avoids this dichotomy. At the algebraic boundary, the squared components remain exactly equal to the mean component weight $1/(N+1)$, thereby preserving finite normalization without introducing artificial damping factors.

Under the unitary transformation into the position basis, the above-established invariance of the boundary term causes the last term in our action difference $D_N = diag\left(1,1,...,1,-N\right)$ to be distributed uniformly over the entire system, resulting in the homogeneous structure of the operator $D_N^{(q)} = I_N - J_N$. Within the present framework, this operator replaces the classical canonical commutation relation. Its matrix form is characterized by the fact that the diagonal elements are identically zero, while all off-diagonal elements take the constant value of minus one. The disappearance of the index $N$ from the algebraic core components of $D_N^{(q)}$ reveals a dimension-independent algebraic structure, which encourages us to refer to $D^{(q)} = I - J$ as the *Leibniz operator*.

From the structure of the Leibniz operator, the position representation of the momentum matrix is derived in Section 2.10. The resulting momentum matrix possesses a discrete Cauchy kernel with purely imaginary off-diagonal elements. For the diagonal elements, the corresponding matrix-element formula yields an indeterminate expression, which is resolved algebraically to zero by imposing matrix antisymmetry. Unlike the continuous Cauchy kernel, whose diagonal singularity requires an external regularization such as the Cauchy principal value or an infinitesimal complex shift, the finite-dimensional algebraic construction determines the diagonal elements intrinsically. The resulting momentum matrix therefore provides a singularity-free algebraic realization of the Cauchy kernel and serves as the algebraic counterpart of the differential momentum operator in the Schrödinger coordinate representation.

### 3.2 The Derivation of Quantum Randomness

The physical measurement problem associated with a sharp position measurement admits a mathematical interpretation through the invariance theorem proven in Section 2.6. In continuous function spaces, a localized state described by a narrow Gaussian wave packet causes the expansion coefficients for high energy values to decay exponentially toward zero. This exponential damping causes the state to appear as if it were essentially concentrated within a low-energy range. In the finite-dimensional formalism presented here, however, the energy components of the normalized position eigenvectors behave fundamentally differently: for each of them, the square of the highest-energy component remains at the mathematical mean $1/(N+1)$ of the squares of all components of the vector. There is no suppression of the highest-energy component relative to the mean component weight. Within the present framework, a sharp position measurement necessarily involves the full informational content of the position state.

At this juncture, the emergence of statistical behaviour can be understood through a direct structural comparison with the classical mechanics of Ludwig Boltzmann [15] and Josiah Willard Gibbs [16]. In classical thermodynamics, randomness arises from the transition from a high-dimensional microstate in phase space to a low-dimensional macroscopic state, resulting in a loss of deterministic information. This classical randomness is epistemic: the underlying equations of motion are strictly deterministic, and the statistical character arises merely from the practical impossibility of recording the immense number of microscopic variables (on the order of $10^{23}$) . In our algebraic formalism, the underlying mechanism of information loss is closely analogous. The physical observer is confronted with a position eigenstate consisting of $N+1$ non-vanishing and undamped components.

Yet, by definition, the observer constitutes a finite system whose measuring apparatus and cognitive capacities can register only finitely many discrete units of information, irrespective of how large the Hilbert-space dimension $N+1$ becomes. Thus, an epistemic origin of randomness naturally emerges. Here, however, lies the essential difference from Boltzmann and Gibbs: whereas the large number of microscopic degrees of freedom is fixed in the classical case, the Hilbert-space dimension is finite but arbitrarily large in the present framework. Consequently, as $N$ increases, only a vanishingly small fraction of the position state's information remains physically accessible to the finite observer. Unlike the classical situation, this informational limitation cannot be overcome by increasing observational resources. We define this irreducible collapse of accessible information as *intrinsic decoherence*. Unlike environmental decoherence, which scatters phase information into an external reservoir, intrinsic decoherence refers to the irreducible loss of accessible information caused by the finite observational window. The resulting indeterminacy constitutes a Boltzmann-like form of epistemic randomness rendered fundamentally irreducible by the arbitrarily large Hilbert-space dimension. Within the present interpretation, quantum randomness is therefore understood not as an intrinsic lawlessness of nature, but as the unavoidable statistical manifestation of restricted information accessibility during a sharp position measurement.

In the continuous theory, this intrinsic decoherence remains concealed during position measurements: owing to the exponential decay of the components, the finite observer can effectively identify the bounded observational domain with the low-energy sector, thereby masking the underlying information loss. However, there is an even deeper reason why randomness is fundamentally unrecognizable within von Neumann's Hilbert-space framework. Although each finite-dimensional space $\mathbb{C}^{N+1}$ embeds into the separable Hilbert space, its distinguished boundary structure is absorbed by Hilbert-space completion and no longer exists as an explicit algebraic object. The finite observational window therefore disappears from the formalism, and the standard formulation introduces quantum randomness through the independent axioms of the projection postulate and the Born probability rule.

The informational collapse within the position measurement stands in sharp contrast to a measurement of energy within the native Fock basis. The determination of the discrete energy eigenvalues proceeds entirely without such difficulties. Since the energy basis represents the native countable structure of the algebra, no projection losses over arbitrarily large spaces occur here. A finite observer can register the discrete energy states within the bulk error-free and deterministically, without suffering a loss of accessible information. Within the present framework, randomness is therefore not interpreted as a universal property of reality, but as a phenomenon of the spatial representation arising from the limited capabilities of the finite observer during a sharp position measurement.

### 3.3 The Leibniz Operator and the Emergence of Relational Space

This section examines the system's underlying relational structure through the Leibniz operator. In the position representation, the operator is characterized by vanishing diagonal entries and uniform off-diagonal entries equal to minus one. Within the present framework, this algebraic structure admits a natural physical interpretation: the vanishing diagonal entries indicate that isolated position values possess no intrinsic significance with respect to the action difference, whereas the uniform off-diagonal entries describe an equal coupling between every pair of position states, independent of their spatial separation. Accordingly, the present framework no longer treats space as an absolute Newtonian stage. Instead, spatial structure emerges from the network of relations among position states, providing a mathematical realization of Leibniz's relational conception of space.

Because the commutator associated with the Leibniz operator scales linearly with $\hbar$, it vanishes in the classical limit as $\hbar$ approaches zero. In this limit, the algebraic coupling encoded by the commutator disappears, allowing the classical description of space as Newton's absolute stage to be recovered.

## 4. OUTLOOK

For large dimensions, the zeros of Hermite polynomials become approximately equidistant in the central spectral region, supporting the lattice model of Appendix C. A natural direction for further research is to investigate the corresponding relation between the finite-dimensional momentum matrices and the Schrödinger momentum operator for the actual nonuniform spectra.

The present work raises questions concerning the origin of quantum randomness. Within the proposed algebraic framework, the limited information accessible to a finite observer suggests a possible epistemic interpretation. Further investigation is needed to clarify the implications of this perspective.

## Appendix A: Inductive Proof of the Eigenvector Component Formula

Derivation of the unnormalized eigenvector components $\bar{u}_{l,k}^{(N)}(\lambda)$ of the tridiagonal position $(N+1)\times(N+1)$ matrix $Q_N$, introduced in Section 2.1.

In this chapter, we omit the superscript $(N)$ for brevity.

**Theorem**

With the initial boundary condition $\bar{u}_{l,0}=1$, the components of any unnormalized eigenvector

$\bar{u}_l$ with $l=0,1,...,N$ are given explicitly by: $$\bar{u}_{l,k}(\lambda_l)=\frac{H_k(\lambda_l)}{\sqrt{2^k\cdot k!}} \tag{A.1}$$

where $H_k(\lambda)$ denotes the $k$-th classic Hermite polynomial, with the recurrence relation $H_{k+1}(\lambda)=2\lambda\cdot H_k(\lambda)-2k\cdot H_{k-1}(\lambda)$.

**Proof by Induction**

The matrix equation $Q_N \bar{u}_l = \lambda_l \bar{u}_l$ dictates the following three-term recurrence relation for the components: $\sqrt{k+1} \cdot \bar{u}_{l,k+1} + \sqrt{k} \cdot \bar{u}_{l,k-1} = \sqrt{2} \cdot \lambda_l \cdot \bar{u}_{l,k}$

Solving this relation for the next component yields the inductive step formula:

$$\bar{u}_{l,k+1} = \frac{\sqrt{2}\,\lambda_l \cdot \bar{u}_{l,k} - \sqrt{k} \cdot \bar{u}_{l,k-1}}{\sqrt{k+1}} \tag{A.2}$$

**Base Case** ($\boldsymbol{k = 0}$ and $\boldsymbol{k = 1;\ l = 0,1,\dots,N}$)**:**

For $k = 0$: Since $H_0(\lambda) = 1$, the formula yields $\bar{u}_{l,0} = \frac{1}{\sqrt{2^0 \cdot 0!}} = 1$.

This satisfies the initial condition.

For $k = 1$ : The first row of the matrix equation gives $\bar{u}_{l,1} = \sqrt{2}\,\lambda_l \cdot \bar{u}_{l,0} = \sqrt{2}\,\lambda_l$ Evaluating our formula for $k = 1$ with $H_1(\lambda) = 2\lambda$ yields: $\bar{u}_{l,1} = \frac{2\lambda_l}{\sqrt{2^1 \cdot 1!}} = \sqrt{2}\,\lambda_l$ . Thus, the base cases hold true.

**Inductive Step** ($\boldsymbol{k \to k+1}$)**:**

We assume that the formula holds for the components $\bar{u}_{l,k}$ and $\bar{u}_{l,k-1}$ (Inductive Hypothesis). Substituting these expressions into the right-hand side of Eq. (A.1) gives:

$$\bar{u}_{l,k+1} = \frac{1}{\sqrt{k+1}} \left( \frac{\sqrt{2}\,\lambda_l \cdot H_k(\lambda_l)}{\sqrt{2^k \cdot k!}} - \frac{\sqrt{k} \cdot H_{k-1}(\lambda_l)}{\sqrt{2^{k-1} \cdot (k-1)!}} \right) = \frac{\sqrt{2}\,\lambda_l \cdot H_k(\lambda_l)}{\sqrt{2^k \cdot k! \cdot (k+1)}} - \frac{\sqrt{k} \cdot H_{k-1}(\lambda_l)}{\sqrt{2^{k-1} \cdot (k-1)! \cdot (k+1)}}$$

Multiplying the first numerator by $\sqrt{2}$ and the second term by $2\sqrt{k}$ yield $\bar{u}_{l,k+1} = \frac{2\lambda_l \cdot H_k(\lambda_l) - 2k \cdot H_{k-1}(\lambda_l)}{\sqrt{2^{k+1} \cdot (k+1)!}} = \frac{H_{k+1}(\lambda_l)}{\sqrt{2^{k+1} \cdot (k+1)!}}$. This completes the inductive step.

By mathematical induction, the formula holds for all $k = 0,1,\dots,N$. This completes the proof.

**Appendix B**

For the classical Hermite polynomials $H_j(x)$ the general Christoffel-Darboux identity for the sum of weighted squares up to order N is given by:

$$\sum_{j=0}^{N} \frac{H_j(x)^2}{2^j\, j!} = \frac{1}{2^{N+1} N!} \left( H'_{N+1}(x) H_N(x) - H_{N+1}(x) H'_N(x) \right).$$

In the present model, the arguments $x = \lambda_l$ are the eigenvalues of the position matrix, which mathematically correspond exactly to the roots of the next-order Hermite polynomial: $H_{N+1}(\lambda_l) = 0, \forall l = 0,\dots,N$ .

Substituting this condition into the right-hand side of the general identity causes the second term inside the bracket to vanish:

$$\sum_{j=0}^{N}\frac{H_j(\lambda_l)^2}{2^j j!}=\frac{1}{2^{N+1}N!}H'_{N+1}(\lambda_l)H_N(\lambda_l).$$

To transform the derivative $H'_{N+1}(\lambda_l)$ we employ the standard recurrence relation for the derivatives of Hermite polynomials: $H'_n(x)=2n\cdot H_{n-1}(x)$.
Setting the index to $n=N+1$ yields: $H'_{N+1}(\lambda_l)=2(N+1)H_N(\lambda_l)$.

Substituting this expression for the derivative back into our simplified identity gives:

$$\sum_{j=0}^{N}\frac{H_j(\lambda_l)^2}{2^j j!}=\frac{1}{2^{N+1}N!}\big(2(N+1)H_N(\lambda_l)\big)\cdot H_N(\lambda_l)=\frac{(N+1)\cdot H_N(\lambda_l)^2}{2^N\cdot N!}$$

Thus all $N$-components are $u_{l,N}^{(N)}=\dfrac{\dfrac{H_N(\lambda_l)}{\sqrt{2^N\cdot N!}}}{\sqrt{\sum_{j=0}^{N}\dfrac{H_j(\lambda_j)^2}{2^j\cdot j!}}}=\dfrac{\dfrac{H_N(\lambda_l)}{\sqrt{2^N\cdot N!}}}{\sqrt{\dfrac{(N+1)\cdot H_N(\lambda_l)^2}{2^N\cdot N!}}}=\dfrac{1}{\sqrt{N+1}}$.

**Appendix C: Relation to the Schrödinger Representation**

The Wintner–Wielandt theorem provides a natural point of departure for two distinct approaches to quantum mechanics. The conventional Schrödinger representation employs unbounded operators in an infinite-dimensional Hilbert space, whereas the second route pursued in the present work retains finite, albeit arbitrarily large, matrices at the price of modifying the canonical commutation relation. The finite-dimensional construction establishes a sequence of algebraic relations: the modified canonical commutation relation gives rise to the action difference, whose Leibniz form in the position representation determines the Cauchy structure of the momentum matrix. The purpose of the following calculation is to show how this structure connects to the Schrödinger momentum operator. Diagonalization of the position matrix $Q_N$ gives $Q_N^{(q)}=diag(q_0,q_1,\dots,q_N)$ so that $\left(Q_N^{(q)}\psi\right)_j=q_j\psi_j$. This is the finite-dimensional counterpart of the Schrödinger multiplication operator $(\hat{q}\psi)(q)=q\psi(q)$. In the position basis, the transformed commutator satisfies $\left[Q_N^{(q)},P_N^{(q)}\right]/i=D_N^{(q)}=I_N-J_N$, where $J_N$ denotes the $(N+1)\times(N+1)$ matrix whose entries are all equal to unity. The appearance of $J_N$ is a direct consequence of the upper boundary term. As derived in Eqs. (13) and (14), the momentum operator in the position basis has the Cauchy form $(P_N^{(q)})_{jk}=\dfrac{i}{q_k-q_j},\ \ j\neq k$ with $(P_N^{(q)})=0$. The validity of this form relies on the presence of $J_N$ in the transformed commutator, which originates from the finite-dimensional boundary term. To establish a more explicit connection with the Schrödinger momentum operator, we introduce the diagonal unitary matrix $\Phi_N=diag(1,-1,1,-1,\dots,(-1)^N)$, which represents an alternating phase transformation. We transform both position $\tilde{Q}_n^{(q)}=\Phi_N^{\dagger}Q_N^{(q)}\Phi_N$ and momentum $\tilde{P}_n^{(q)}=\Phi_N^{\dagger}P_N^{(q)}\Phi_N$. Since $Q_N^{(q)}$ and $\Phi_N$ are diagonal, the position matrix remains unchanged: $\tilde{Q}_n^{(q)}=Q_N^{(q)}=diag(q_0,q_1,\dots,q_N)$.

The momentum matrix, however, acquires alternating signs: $\left(\tilde{P}_N^{(q)}\right)_{jk} = \dfrac{i\cdot(-1)^{j+k}}{q_k - q_j}$, $j \neq k$, while $\left(\tilde{P}_N^{(q)}\right)_{jj} = 0$. Since $(-1)^{j+k} = (-1)^{k-j}$, the signs of the entries on odd-numbered off-diagonals are reversed, whereas those of the even off-diagonals remain unchanged. We now apply the transformed momentum matrix to a state vector $\psi = (\psi_0, \psi_1, \ldots, \psi_N)$ : $\left(\tilde{P}_N^{(q)}\psi\right)_j = i \sum_{k=0,k\neq j}^{N} \dfrac{(-1)^{j+k}}{q_k - q_j}\psi_k$. As $N$ increases, the zeros of the Hermite polynomials become increasingly equidistant over a central region of the position spectrum whose width can grow without bound, motivating the equidistant lattice model considered below. To elucidate the connection with the continuum representation, we consider the simplifying approximation of equally spaced position eigenvalues $q_j = j\,\Delta q$ , $q_k = k\,\Delta q$ . The equation then becomes $\left(\tilde{P}_N^{(q)}\psi\right)_j = \dfrac{i}{\Delta q}\sum_{k=0,k\neq j}^{N} \dfrac{(-1)^{j+k}}{k-j}\psi_k$. Grouping the terms with $k = j + m$ and $k = j - m$, where $m \geq 1$ , we combine the contributions from points symmetrically located about $q_j$. Extending the equidistant lattice to infinity, we obtain $\left(\tilde{P}^{(q)}\psi\right)_j = \dfrac{i}{\Delta q}\sum_{m\geq 1}\dfrac{(-1)^m}{m}(\psi_{j+m} - \psi_{j-m})$. Consider first the contribution from the nearest neighbours, $m = 1$: $\left(\tilde{P}^{(q)}\psi\right)_j^{(m=1)} = -\dfrac{i}{\Delta q}(\psi_{j+1} - \psi_{j-1})$. For a sufficiently smooth wave function and small $\Delta q$ , the central difference quotient approaches the derivative, $\dfrac{\psi_{j+1} - \psi_{j-1}}{2\,\Delta q} \approx \dfrac{d\psi}{dq}(q_j)$. Consequently, $\left(\tilde{P}^{(q)}\psi\right)_j^{(m=1)} \approx -2i\dfrac{d\psi}{dq}(q_j)$. This result has the correct sign but is twice the Schrödinger momentum operator. The remaining off-diagonal contributions are therefore essential for obtaining the correct coefficient. To evaluate the complete sum, we consider a plane wave $\psi(q) = e^{ipq}$ where $p$ denotes wave number. On such lattice, we have $\psi_{j+m} = e^{i(q_j + m\Delta q)} = e^{ipm\Delta q}\psi_j$ and $\psi_{j-m} = e^{ip(q_j - m\Delta q)} = e^{-ipm\Delta q}\psi_j$.

Their difference is $\psi_{j+m} - \psi_{j-m} = \left(e^{ipm\Delta q} - e^{-ipm\Delta q}\right)\psi_j = 2i\sin(mp\Delta q)\psi_j$ . Substitution into the above equation yields $\left(\tilde{P}^{(q)}\psi\right)_j = -\dfrac{2}{\Delta q}\sum_{m\geq 1}\dfrac{(-1)^m}{m}\sin(mp\Delta q)\psi_j$ . The Fourier series identity $\sum_{m=1}^{\infty}\dfrac{(-1)^m}{m}\sin(mx) = -\dfrac{x}{2}$ , $-\pi < x < \pi$ allows the infinite sum to be evaluated explicitly. Setting $x = p\Delta q$ , we obtain $\left(\tilde{P}^{(q)}\psi\right)_j = -\dfrac{2}{\Delta q}\left(-\dfrac{p\Delta q}{2}\right)\psi_j = p\psi_j$ , provided that $|p\Delta q| < \pi$. This is precisely the action of the Schrödinger momentum operator on a plane wave: $\hat{p}\psi(q) = -i\dfrac{d}{dq}e^{ipq} = pe^{ipq} = p\psi(q)$.

Thus, on the idealized infinite equidistant lattice, the alternating Cauchy matrix reproduces exactly the action of the Schrödinger momentum operator on plane waves for $\left|p\,\Delta q\right| < \pi$. This establishes an explicit structural connection between the finite-dimensional algebraic construction and the analytic representation.

On the idealized infinite equidistant lattice, the resulting alternating Cauchy matrix coincides, up to the factor $-i$ , with the sinc differentiation matrix known from spectral methods [17]. Here, this structure emerges from the finite-dimensional algebraic construction through the action difference and the alternating phase transformation.

**Acknowledgements** We would like to thank Gundula Fromm, physician and psychologist, whose questions and suggestions have repeatedly inspired and challenged the preparation of this work.